\documentclass[12pt]{article}
\usepackage{amssymb, amsfonts, amsmath, amsthm, url, graphicx}

\def\3{\subset }
\def\4{\subseteq }

\def\calp{{\cal P}}
\def\0{\leqno}

\def\barr{\begin{array}}
\def\earr{\end{array}}
\def\dd{\displaystyle}

\def\Z{{\rlap{$\kern2pt{\rm Z}$}{\rm Z}\,}}

\title{A characterization of\\ elementary abelian 2-groups}
\author{Marius T\u arn\u auceanu}
\date{November 27, 2016}

\begin{document}

\maketitle

\begin{abstract}
In this note we give a characterization of elementary
abelian 2-groups in terms of their maximal sum-free subsets.
\end{abstract}

\noindent{\bf MSC (2010):} Primary 11B75; Secondary 20D60, 20K01.

\noindent{\bf Key words:} elementary abelian 2-groups, sum-free subsets,
maximal sum-free subsets, ma\-xi\-mal subgroups.

\section{Introduction}

Throughout this paper $G$ will denote an additive, but not necessarily abelian, group.
A subset $A$ of $G$ is said to be \textit{sum-free} if $(A+A)\cap A=\emptyset$, that is the equation $x+y=z$ has no solutions $x,y,z\in A$. Sum-free sets were introduced by Schur \cite{12} whose celebrated result states that the set of positive integers cannot be partitioned into finitely many sum-free sets. Another famous result on sum-free sets is a conjecture of Cameron and Erd\H{o}s \cite{3,4}, which asserts that the number of sum-free subsets of $\{1,2,...,n\}$ is $O(2^{n/2})$. It has been recently proved by Green \cite{7}.

Sum-free subsets of abelian groups are central objects of interest in Additive Combinatorics, and have been studied intensively in the last years (see e.g. \cite{1,2}, \cite{5,6} and \cite{8}-\cite{11}). Probably the most important problems in this direction are: \textit{How large a sum-free subset of $G$ can be}\,? and \textit{How many sum-free sets of $G$ are there}\,?

In the following we will focus on maximal sum-free subsets of $G$. Our main theorem characterizes elementary abelian $2$-groups by connecting these subsets with the complements of maximal subgroups.
\newpage

\bigskip\noindent{\bf Theorem 1.1.} {\it A finite group $G$ is an elementary abelian $2$-group
if and only if the set of maximal sum-free subsets coincides with the set of complements
of maximal subgroups.}
\bigskip

The following corollary is an immediate consequence of Theorem 1.1.

\bigskip\noindent{\bf Corollary 1.2.} {\it The number of maximal sum-free subsets of \,$\mathbb{Z}_2^n$ is $2^n-1$,
and each of them has $2^{n-1}$ elements.}
\bigskip

Theorem 1.1 can be also used to count the sum-free subsets of $\mathbb{Z}_2^n$ since these are all subsets of maximal sum-free subsets. Let $k=2^n-1$ and denote by $M_1, M_2, ..., M_k$ the maximal subgroups of $\mathbb{Z}_2^n$. Then, by applying the well-known Inclusion-Exclusion Principle, we infer that the total number of sum-free subsets of $\mathbb{Z}_2^n$ is
$$x_n=|\bigcup_{i=1}^k \calp(\mathbb{Z}_2^n\setminus M_i)|=\dd\sum_{r=1}^k (-1)^{r-1}\hspace{-7mm}\dd\sum_{1\leq i_1<\cdots<i_r\leq k}\hspace{-5mm}|\calp(\mathbb{Z}_2^n\setminus M_{i_1})\cap\cdots\cap\calp(\mathbb{Z}_2^n\setminus M_{i_r})|$$
$$=\dd\sum_{r=1}^k (-1)^{r-1}\hspace{-7mm}\dd\sum_{1\leq i_1<\cdots<i_r\leq k}\hspace{-5mm}|\calp(\mathbb{Z}_2^n\setminus M_{i_1}\cup\cdots\cup M_{i_r})|=\dd\sum_{r=1}^k (-1)^{r-1}\hspace{-7mm}\dd\sum_{1\leq i_1<\cdots<i_r\leq k}\hspace{-6mm}2^{2^n-|M_{i_1}\cup\,\cdots\,\cup M_{i_r}|}.$$
\smallskip

We remark that this depends only on the cardinality of any union of $M_i$'s. Note that its first values are: $x_1=2$, $x_2=7$, $x_3=64$, $x_4=3049$, ... and so on. Thus, a natural problem is the following.

\bigskip\noindent{\bf Open problem.} Give an explicit formula for $x_n$.
\bigskip

Most of our notation is standard and will usually not be repeated
here. Given a subset $A$ of $G$, we write $A+A=\{a_1+a_2\,|\,a_1,a_2\in A\}$, $A-A=\{a_1-a_2\,|\,a_1,a_2\in A\}$, $\dd\frac{1}{2}A=\{a\in G\,|\,2a\in A\}$ and $G\setminus A=\{x\in G\,|\, x\notin A\}$.

\section{Proofs of the main results}

We start with two easy but important lemmas.

\bigskip\noindent{\bf Lemma 2.1.} {\it Let $G$ be a group and $H$ be a subgroup of $G$. Then
$G\setminus H$ is a non-empty maximal sum-free subset of $G$ if and only if $[G:H]=2$.}

\bigskip\noindent{\bf Proof.} First of all, we remark that if $A$ is a sum-free of $G$, then $|A|\leq\dd\frac{|G|}{2}$\,. Indeed, given $x\in A$ the set $x+A=\{x+y\,|\,y\in A\}$ has the same number of elements as $A$ and $(x+A)\cap A=\emptyset$. In other words, $G$ contains at least $2|A|$ elements.

Assume first that $G\setminus H$ is a non-empty maximal sum-free subset of $G$. By the above remark we have $|G\setminus H|\leq\dd\frac{|G|}{2}$\,, or equivalently $[G:H]\leq 2$. This leads to $[G:H]=2$ because $G\setminus H$ is non-empty.

Conversely, assume that $[G:H]=2$ and let $x,y\in G\setminus H$. Since $x+H$ and $-y+H$ are non-trivial cosets, it
follows that $x+H=-y+H$, that is $x+y\in H$. Therefore $G\setminus H$ is a sum-free subset of $G$. Now, let $x\in G\setminus H$ and $z\in H$. Then $x+z\in G\setminus H$. This shows that $(G\setminus H)\cup\{z\}$ is not a sum-free subset of $G$, and consequently $G\setminus H$ is a non-empty maximal sum-free subset of $G$, as desired.
\hfill\rule{1,5mm}{1,5mm}

\bigskip\noindent{\bf Lemma 2.2.} {\it Let $G$ be a group and $A$ be a subset of $G$. Then
$A$ is a maximal sum-free subset of $G$ if and only if $G\setminus A=(A+A)\cup (A-A)\cup\dd\frac{1}{2}A$.}

\bigskip\noindent{\bf Proof.} Assume first that $A$ is a maximal sum-free subset of $G$. Then it is a routine to check that $(A+A)\cup (A-A)\cup\dd\frac{1}{2}A\subseteq G\setminus A$. Let $x\in G\setminus A$. By our hypothesis, we infer that $A\cup\{x\}$ is not a maximal sum-free subset of $G$ and so there are $a,b\in A\cup\{x\}$ such that $a+b\in A\cup\{x\}$. If $a=b=x$ we have either $x+x\in A$, i.e. $x\in\dd\frac{1}{2}A$, or $x+x=x$, i.e. $x=0\in A-A$. If $a\in A$ and $b=x$, then $a+x\neq x$ and so $a+x\in A$, i.e. $x\in A-A$. Note that a similar conclusion follows for $b\in A$ and $a=x$. Finally, if $a,b\in A$, then we must have $a+b=x$, i.e. $x\in A+A$. These show that $G\setminus A\subseteq(A+A)\cup (A-A)\cup\dd\frac{1}{2}A$.

Conversely, assume that $G\setminus A=(A+A)\cup (A-A)\cup\dd\frac{1}{2}A$. Since $A+A\subseteq G\setminus A$, we infer that $A$ is a sum-free subset of $G$. It remains to prove that this is maximal, that is $A\cup\{x\}$ is not a sum-free subset of $G$ for any $x\in G\setminus A$. By our hypothesis, we distinguish the following three cases. If $x\in A+A$, then $x=a+b$ for some $a,b\in A$; thus $a,b\in A\cup\{x\}$ and $a+b\in A\cup\{x\}$. If $x\in A-A$, then $x=a-b$ for some $a,b\in A$; thus $x,b\in A\cup\{x\}$ and $x+b=a\in A\cup\{x\}$. Finally, if $x\in\dd\frac{1}{2}A$, then $x\in A\cup\{x\}$ and $x+x\in A\subseteq A\cup\{x\}$. This completes the proof.
\hfill\rule{1,5mm}{1,5mm}
\bigskip

We are now able to prove Theorem 1.1.

\bigskip\noindent{\bf Proof of Theorem 1.1.}
Let $G=\mathbb{Z}_2^n$ for some $n\geq 1$. Then all maximal subgroups of $G$ are of index $2$, and therefore their complements are maximal sum-free subsets by Lemma 2.1. Suppose now that $A$ is a maximal sum-free subset of $G$. Then $x+x=0\notin A \,\forall\, x\in G$, that is $\dd\frac{1}{2}A=\emptyset$, and $A+A=A-A$. These lead to $$G\setminus A=A+A$$in view of Lemma 2.2. Let $\{e_1,e_2,...,e_m\}$ be a maximal linearly independent subset of $A$ (over $\mathbb{Z}_2$). We infer that $A$ consists of all sums of an odd number of $e_i$'s. If $m<n$ we can choose $e_{m+1}\in G$ such that the vectors $e_1,e_2,...,e_{m+1}$ are linearly independent. One obtains that
$$A'=\{e_1,e_2,...,e_{m+1},e_1+e_2+e_3,e_1+e_2+e_4,...,e_{m-1}+e_m+e_{m+1},...\}$$is a sum-free subset of $G$ and $A\subsetneq A'$, contradicting the maximality of $A$. Thus $m=n$ and
$$|A|=\binom{n}{1}+\binom{n}{3}+\cdots+\binom{n}{2r+1}+\cdots=2^{n-1},$$which implies that
$$|G\setminus A|=|A+A|=2^{n-1}.$$On the other hand, $G\setminus A$ is a subgroup of $G$ because it consists of all sums of an even number of $e_i$'s, and so it is a maximal subgroup.

Conversely, let $M_1, M_2, ..., M_k$ be the maximal subgroups of $G$ and suppose that $G\setminus M_i$, $i=1,2,...,k$, are the maximal sum-free subsets of $G$. Then $[G:M_i]=2$, for any $i=1,2,...,k$, by Lemma 2.1. We infer that $G$ is a nilpotent group, more precisely a $2$-group. Since every non-trivial element of $G$ is contained in a maximal sum-free subset of $G$, we have
$$G\setminus 1=\bigcup_{i=1}^k\, G\setminus M_i=G\setminus\,\bigcap_{i=1}^k M_i=G\setminus\Phi(G),$$that is $\Phi(G)=1$. Hence $G$ is an elementary abelian $2$-group.
\hfill\rule{1,5mm}{1,5mm}

\bigskip\noindent{\bf Acknowledgements.} The author is grateful to the reviewer for
its remarks which improve the previous version of the paper.

\vspace*{5ex}\small

\hfill
\begin{minipage}[t]{5cm}
Marius T\u arn\u auceanu \\
Faculty of  Mathematics \\
``Al.I. Cuza'' University \\
Ia\c si, Romania \\
e-mail: {\tt tarnauc@uaic.ro}
\end{minipage}

\end{document}